\documentclass[11pt,a4paper]{amsart}
\usepackage{hyphenat}
\usepackage[utf8]{inputenc}
\usepackage{lmodern}

\usepackage[T1]{fontenc}
\usepackage{microtype}
\usepackage[normalem]{ulem}
\usepackage{amsthm}
\usepackage{amsmath,amssymb,amsfonts,amsthm,color,latexsym}
\usepackage{xcolor}
\usepackage{float}
\usepackage{xcolor}
\usepackage{comment}

\providecommand{\cprime}{'}
\usepackage{booktabs}
\newcommand{\C}{\mathbb{C}}

\newtheorem{theorem}{Theorem}[section]

\usepackage[colorlinks=true,pagebackref=true]{hyperref}
\hypersetup{
  linkcolor=blue,   
  citecolor=blue,   
  urlcolor=blue     
}

\newcommand{\cpt}[1]{\mathbb{C}P^{2}}

\newcommand{\gsv}{\mbox{\rm GSV}}

\newtheorem{lemma}[theorem]{Lemma}
\newtheorem{proposition}[theorem]{Proposition}
\newtheorem{corollary}[theorem]{Corollary}
\newtheorem{conjecture}[theorem]{Conjecture}

\newtheorem{example}[theorem]{Example}
\newtheorem{remark}[theorem]{Remark}

\title[Higher-order Milnor and Tjurina numbers]{ The $k$-th Milnor and Tjurina numbers: conjectures and a bifurcation formula for pairs}
\date{\today}

\author[A. Fern\'andez-P\'erez]{Arturo Fern\'andez-P\'erez}
\address[Arturo Fern\'{a}ndez P\'erez] {Department of Mathematics. Federal University of Minas Gerais. Av. Pres. Ant\^onio Carlos, 6627 
CEP 31270-901\\
Pampulha - Belo Horizonte - Brazil. ORCID: 0000-0002-5827-8828}
\email{fernandez@ufmg.br}

\author[E. R.  Garc\'{i}a Barroso]{Evelia R. Garc\'{i}a Barroso}
\address[Evelia R. Garc\'{i}a Barroso]{Departamento de Matem\'{a}ticas, Estad\'{\i}stica e Investigaci\'on Ope\-rativa\\
Instituto Universitario de Matem\'aticas y Aplicaciones (IMAULL). Universidad de La Laguna. Apartado de Correos 456. 38200 La Laguna, Te\-ne\-rife, Spain. ORCID: 0000-0001-7575-2619 }
\email{ergarcia@ull.es}

\author[M.~Ribeiro]{Marcela Ribeiro}
\address[Marcela Ribeiro]{Department of General Education, Universidade Federal de Ciência e Inovação de Minas Gerais – Campus Curvelo.
Rua Raymundo Mattoso, 900
Santa Rita - Curvelo - MG 35790-636\\
Brazil. ORCID: 0009-0009-2047-6341}
\email{marcelaribeiro@cefetmg.br}

\author[N. Saravia-Molina]{Nancy Saravia-Molina}
\address[Nancy Saravia-Molina]{Secci\'{o}n de Matem\'{a}ticas, Departamento Acad\'emico de Ciencias, Pontificia Universidad Cat\'{o}lica del Per\'{u}, Av. Universitaria 1801,
San Miguel, Lima 32, Peru. ORCID: 0000-0002-2819-8835}
\email{nsaraviam@pucp.edu.pe}

\thanks{The first author acknowledges
support from CNPq Projeto Universal
408687/2023-1 ``Geometria das
Equa\c{c}\~oes Diferenciais
Alg\'ebricas'' and from CNPq-Brazil
PQ-306011/2023-9.
The second author is supported by the
Spanish grant PID2023-149508NB-I00,
funded by MICIU/AEI/10.13039/501100011033
and by FEDER, UE.
The third author acknowledges support
from CAPES, Brazil.
The first and fourth authors were also
supported by the Vicerrectorado de
Investigaci\'on (VRI) of the Pontificia
Universidad Cat\'olica del Per\'u (PUCP)
through grant DFI-2025-PI1275.}


\usepackage{lineno}

\usepackage{todonotes}

\begin{document}

\begin{abstract}
We study the higher-order Milnor and
Tjurina numbers of isolated hypersurface
singularities and extend these invariants
to pairs of plane holomorphic germs.
We first show that a conjectured upper
bound of Hussain, Yau and Zuo on the ratio
between the higher-order Milnor and
Tjurina numbers of isolated hypersurface
singularities in three-dimensional complex
space does not hold: it fails already for
the simplest singularity, a non-degenerate
quadratic one, and the ratio tends to
infinity as the order grows.
We then prove that a conjecture of Shi,
Yau and Zuo on semiquasihomogeneous
singularities holds for plane curves, for
every order and over formal power series,
and show that it fails in dimension three.
Using the theory of higher-order Milnor
and Tjurina numbers of foliations in
dimension two, we introduce these
invariants for pairs of plane germs.
For every order, we show that the
higher-order Milnor number of a pair is
determined by the Milnor number of the
pair and the algebraic multiplicity of
the associated foliation, and we compute
the difference between the higher-order
Milnor and Tjurina numbers of the pair.
Finally, we establish a bifurcation
formula comparing the higher-order Milnor
number of a pair with that of the product
of the germs. As a consequence, we
characterize semitame meromorphic germs as
those for which these two numbers agree at
some order, or equivalently at every order.
\end{abstract}

\maketitle

\tableofcontents
\section{Introduction}

Let $\mathcal{O}_n=\mathbb{C}\{x_1,\ldots,x_n\}$ denote
the ring of germs of holomorphic functions at the origin of
$\mathbb{C}^n$, with maximal ideal $\mathfrak{m}$. Let
$f\in\mathfrak{m}$ define an isolated hypersurface singularity
$(V,0)=(f^{-1}(0),0)$, and denote by
$J(f)=(f_{x_1},\ldots,f_{x_n})$ its Jacobian ideal.
The Milnor and Tjurina numbers of $f$ are defined by
$\mu(f)=\dim_{\mathbb{C}}\mathcal{O}_n/J(f)$ and
$\tau(f)=\dim_{\mathbb{C}}\mathcal{O}_n/((f)+J(f))$,
respectively.
The Milnor number is the rank of the middle homology of the
Milnor fiber of $f$ \cite{Milnor68}, whereas the Tjurina number
is the dimension of the base of a semi-universal deformation
of $(V,0)$ \cite{Tjurina69}. One always has
$\tau(f)\leq\mu(f)$, and Saito's theorem \cite{Saito} states
that equality holds if and only if $f$ is quasihomogeneous
up to a change of coordinates. Moreover, the Mather--Yau
theorem \cite{MatherYau} asserts that the algebra
$\mathcal{O}_n/((f)+J(f))$ determines the analytic type
of $(V,0)$.

The algebras
$\mathcal{O}_n/\mathfrak{m}^kJ(f)$ and
$\mathcal{O}_n/((f)+\mathfrak{m}^kJ(f))$
were first considered by Greuel--Pham \cite{Pham}.
The corresponding $k$-th Milnor and Tjurina numbers
were subsequently studied by Hussain--Liu--Yau--Zuo
\cite{k-Milnor} and Hu--Yau--Zuo \cite{Hu-Yau-Zuo}.
For each integer $k\geq 0$, one considers
$M_k(f)=\mathcal{O}_n/\mathfrak{m}^kJ(f)$ and
$T_k(f)=\mathcal{O}_n/((f)+\mathfrak{m}^kJ(f))$.
Their dimensions
$\mu^k(f)=\dim_{\mathbb{C}}M_k(f)$ and
$\tau^k(f)=\dim_{\mathbb{C}}T_k(f)$
are called the $k$-th Milnor number and the $k$-th Tjurina
number of $f$, respectively. For $k=0$, these are the
classical Milnor and Tjurina numbers, and
$\tau^k(f)\leq\mu^k(f)$ for every $k\geq 0$.

The relation between these numbers differs
from the classical case. If $f$ is quasihomogeneous,
the Euler identity gives $f\in\mathfrak{m}J(f)$, so
$\tau^1(f)=\mu^1(f)$. For $k\geq 2$, however,
$f$ need not belong to $\mathfrak{m}^kJ(f)$, and the
two numbers may differ even for a quasihomogeneous
singularity. Questions concerning the quotient $\mu^k(f)/\tau^k(f)$ and
the relation between the $k$-th Milnor number of a
semiquasihomogeneous germ and that of its quasihomogeneous
principal part have led to several conjectures
\cite{Hussain-Yau,Shi}.

In this paper, we address these questions and investigate
the corresponding invariants for pairs of plane holomorphic
germs. Among our main contributions, we provide counterexamples
to two conjectures concerning isolated hypersurface singularities,
while establishing the validity of one of them for plane curve
singularities. We also derive explicit formulas for the
$k$-th Milnor and Tjurina numbers of pairs and explore
their relationship with the bifurcation fibers of a pencil
of plane curves. Our approach relies on the extension of these invariants to germs of
holomorphic foliations in dimension two, developed by
Ribeiro~\cite{Marcela-tesis} in her doctoral thesis and presented
in~\cite{Marcela}. 
In particular, by applying these results
to the foliation defined by $\omega_{f/g}=gdf-fdg$, we obtain
the desired formulas for pairs. A more detailed study of the Milnor and Tjurina numbers
of foliations on $(\mathbb{C}^2,0)$ can be found in~\cite{Evelia}.

The paper is organized as follows.
Section~\ref{basic} recalls the necessary definitions
and preliminary results.

In Section~\ref{cociente-b}, we investigate a conjecture
of Hussain--Yau--Zuo \cite[Conjecture~1.4]{Hussain-Yau},
which predicts that
$\mu^k(f)/\tau^k(f)<3/2$ for every $k\geq 0$ and
every isolated hypersurface singularity
$V=\{f=0\}\subset(\mathbb{C}^3,0)$.
The case $k=0$ is precisely a conjecture of
Almir\'on~\cite[Conjecture~4.10]{Almiron}, which was
recently proved by Ma--Zuo \cite[Corollary~4.8]{Ma-Zuo}.
We show that the bound does not hold in general for
$k\geq 1$.
Already for $f=x^2+y^2+z^2$ one has
$\mu^k(f)=\binom{k+3}{3}$ and
$\tau^k(f)=(k+1)^2$, so
$\mu^k(f)/\tau^k(f)>3/2$ exactly for $k\geq 5$;
for $k=5$ the quotient is $14/9$.
More generally, for every isolated hypersurface
singularity in $n\geq 2$ variables,
$\mu^k/\tau^k$ tends to infinity as $k\to\infty$;
see Theorem~\ref{thm:growth}.
In particular, the conjecture fails for the families
considered by Hussain--Yau--Zuo, who checked only
$k\leq 2$. For instance, for
$f=x^7+y^7+z^7+x^3y^3z^3$, computations performed
with \textsc{Singular}~\cite{Singular} yield
$\mu^{16}(f)=2024$ and $\tau^{16}(f)=1344$, hence
$\mu^{16}(f)/\tau^{16}(f)=253/168>3/2$.
The analogous conjecture for plane curves
\cite[Conjecture~1.3]{Hussain-Yau} also fails.

Section \ref{k milnor} is devoted to a conjecture
of Shi--Yau--Zuo \cite[Conjecture~4.65]{Shi}
on the behaviour of the $k$-th Milnor numbers
under semiquasihomogeneous perturbations.
More precisely, the conjecture asserts that,
if $f=f_w+f'$ is semiquasihomogeneous with
quasihomogeneous principal part $f_w$, then
$\mu^k(f)=\mu^k(f_w)$ for every $k\geq 0$.
Shi--Yau--Zuo established this equality under
the additional assumption
$k\leq\min_i\{\nu(f_{x_i})\}$
\cite[Theorem~C]{Shi}.
We prove that the conjecture holds for all
$k\geq 0$ in the case of plane curve singularities,
that is, when $f\in\mathbb{C}[[x,y]]$ (and hence also when
$f\in\mathbb{C}\{x,y\}$);
see Theorem~\ref{Teorema-semi}.
Our proof relies on a formula expressing
$\mu^k(f)$ in terms of the classical Milnor number
and the multiplicity of $f$
\cite[Corollary~4.1]{Marcela}, together with the
equality $\mu(f)=\mu(f_w)$, which holds for convergent
series by~\cite{Arnold} and for formal series by
\cite[Proposition~2]{BGM}.
In contrast, we show that the conjecture fails
in dimension three. To this end, we consider
$f=x^3+y^7+z^7+xy^2z^3$ and
$f_w=x^3+y^7+z^7$,
where $f_w$ is quasihomogeneous with weights
$(7,3,3)$. In this case,
$\mu^6(f)=199$ whereas $\mu^6(f_w)=200$,
providing an explicit counterexample. In a recent preprint, Liu--Yau--Zhang
\cite[Main Theorem A]{Liu-Yau-Zhang} proved
the equal-weight case of the conjecture for
convergent germs in every dimension
$n\geq 2$, namely for germs whose
homogeneous principal part defines a smooth
projective hypersurface, equivalently has an
isolated critical point at the origin. Their result complements
Theorem~\ref{Teorema-semi}, which allows arbitrary
positive weights in dimension two. Our
three-variable counterexample has weights $(7,3,3)$
and therefore lies outside their setting.
Their study of higher Milnor sequences also uses
the plane curve formula of
Ribeiro and Fern\'andez-P\'erez
\cite[Corollary~4.1]{Marcela}.

In Section~\ref{k pares}, following
\cite{Tjurina-par,Marcela}, we introduce the $k$-th
Milnor and Tjurina numbers of a pair of plane
holomorphic germs. Let $f,g\in\mathfrak{m}\subset
\mathcal{O}_2$ be reduced and relatively prime, and
assume that the Milnor number of the pair,
defined in \cite{Sz2}, satisfies
\[
\mu(f,g)=
\dim_{\mathbb{C}}
\frac{\mathcal{O}_2}
{(gf_x-fg_x,\;gf_y-fg_y)}
<\infty.
\]
For each integer $k\geq 0$, we define
\[
\mu^k(f,g)=
\dim_{\mathbb{C}}
\frac{\mathcal{O}_2}
{\mathfrak{m}^k(gf_x-fg_x,\;gf_y-fg_y)}.
\]
For a member
$C_{[\alpha:\beta]}:\alpha f+\beta g=0$
of the pencil generated by $f$ and $g$, we also define
\[
\tau^k_{[\alpha:\beta]}(f,g)=
\dim_{\mathbb{C}}
\frac{\mathcal{O}_2}
{\mathfrak{m}^k(gf_x-fg_x,\;gf_y-fg_y)
+(\alpha f+\beta g)}.
\]
Theorem~\ref{teo1} expresses $\mu^k(f,g)$ in terms of
$\mu(f,g)$, $k$, and the algebraic multiplicity of
the $1$-form $\omega_{f/g}=g\,df-f\,dg$.
For every member $C_{[\alpha:\beta]}$ of the
pencil, which is reduced by \cite[Lemma~4.1]{FFMR},
Proposition \ref{teo_2} gives
$\tau^k_{[\alpha:\beta]}(f,g)=
\tau^k(C_{[\alpha:\beta]})+I_0(f,g)$,
where $I_0(f,g)$ denotes the intersection multiplicity
of $f=0$ and $g=0$ at the origin. Thus the difference
between these two Tjurina numbers is independent of
both $k$ and the chosen member of the pencil.
Combining this with the invariant $\eta(f,g)$
introduced by Suwa~\cite{Suwa}, Theorem~\ref{teo3} gives,
for every member $C_{[\alpha:\beta]}$ of the pencil and
every $k\geq 0$,
\[
\mu^k(f,g)-\tau^k_{[\alpha:\beta]}(f,g)=
\eta(f,g)+k(k+1)-\tau^k(C_{[\alpha:\beta]})
-\varepsilon_k(m_0),
\]
where $m_0$ is the algebraic multiplicity of
$\omega_{f/g}$ and $\varepsilon_k(m_0)$ is the correction
term of Theorem~\ref{teo1}.

Finally, in Section~\ref{bifu}, we
compare $\mu^k(f,g)$ with $\mu^k(fg)$.
Starting from the bifurcation formula
\cite[Theorem~4.6]{FFMR}, we express their
difference as the sum of the Milnor number
jumps at the bifurcation fibers distinct from
$f=0$ and $g=0$, together with an explicit
correction depending on $k$ and the
multiplicities of $f$, $g$, and
$\omega_{f/g}$, namely
$\varepsilon_k(M-1)-\varepsilon_k(m_0)$
with $M=\nu(f)+\nu(g)$ and
$m_0=\nu(\omega_{f/g})$; see
Theorem~\ref{kbifu}.
Moreover, $m_0=M-1+\delta$, where
$\delta$ is the sum, over the bifurcation
fibers distinct from $f=0$ and $g=0$, of the
excess of their multiplicity over the generic
one. It vanishes unless $\nu(f)=\nu(g)$ and the
initial forms of $f$ and $g$ are proportional,
and then only one fiber contributes; see
Proposition~\ref{multfol}.
When $\delta=0$, that is, when
$\nu(\omega_{f/g})=\nu(f)+\nu(g)-1$,
the formula reduces to
\[
\mu^k(f,g)-\mu^k(fg)
=
\sum_{[\alpha:\beta]\in\mathcal{B}^*(f,g)}
\bigl(\mu(C_{[\alpha:\beta]})-\mu_{\mathrm{gen}}\bigr)
\qquad\text{for every }k\geq 0.
\]
Here $\mu_{\mathrm{gen}}$ is the Milnor number of
a generic member of the pencil, and
$\mathcal B^*(f,g)$ is the set of bifurcation
parameters different from $[1:0]$ and $[0:1]$.
In this case, this difference is
independent of $k$.
In general, the meromorphic germ
$(f:g)$ is \textit{semitame} if and only if
$\mu^k(f,g)=\mu^k(fg)$ for some, equivalently every integer
$k\geq 0$; see Corollary~\ref{semitame}.
\section{Basic definitions and preliminary results}\label{basic}
Let $(V,0)\subset(\mathbb{C}^n,0)$ be an isolated hypersurface singularity
defined by a germ $f\in\mathcal{O}_n=\mathbb{C}\{x_1,\ldots,x_n\}$, and let
$\mathfrak{m}$ denote the maximal ideal of $\mathcal{O}_n$.
Denote by $\nu(f)$ the order of $f\in\mathcal{O}_n$.
Following Hussain--Liu--Yau--Zuo~\cite{k-Milnor}, for each integer $k\geq 0$
we define the $k$-th Milnor number and the
$k$-th Tjurina number of $f$ by
\[
\mu^k(f)
=\dim_{\mathbb{C}}\frac{\mathcal{O}_n}{\mathfrak{m}^kJ(f)},
\qquad
\tau^k(f)
=\dim_{\mathbb{C}}\frac{\mathcal{O}_n}
{(f)+\mathfrak{m}^kJ(f)},
\]
where $J(f)=(f_{x_1},\ldots,f_{x_n})$ is the Jacobian ideal of
$f$.

In Ribeiro's PhD thesis~\cite{Marcela-tesis} (see also \cite{Marcela}), these notions were
extended to germs of holomorphic foliations in dimension two. More precisely, let
$\omega=P(x,y)\,dx+Q(x,y)\,dy$ be a germ of holomorphic $1$-form defining a
holomorphic foliation $\mathcal{F}$ on $(\mathbb{C}^2,0)$ with an isolated
singularity at the origin, that is, $P,Q\in\mathfrak{m}$ have no common factor.
By definition, the algebraic multiplicity of $\mathcal F$ is
$\nu(\mathcal F):=\min\{\nu(P), \nu(Q)\}$, and we also write $\nu(\omega)=\nu(\mathcal{F})$.

For $f,g\in\mathcal{O}_2$, the intersection multiplicity of $f$ and $g$ at the origin is
$I_0(f,g):=\dim_{\mathbb{C}}\bigl(\mathcal{O}_2/(f,g)\bigr)$,
where $(f,g)$ denotes the ideal generated by $f$ and $g$.
By definition, the Milnor number of $f\in\mathcal{O}_2$ is
$\mu(f)=I_0(f_x,f_y)$, which coincides with $\dim_{\mathbb{C}}\mathcal{O}_2/J(f)$,
and the Milnor number of the foliation $\mathcal{F}$ is
$\mu(\mathcal{F})=I_0(P,Q)$.

Let $C=\{f=0\}$ be a germ of a complex reduced curve, where
$f\in \mathfrak{m}$. The curve $C$ is said to be invariant by
$\mathcal{F}$ or $\mathcal{F}$-invariant if
$(\omega\wedge df)/f\in\Omega^2(\mathbb{C}^2,0)$,
where $\Omega^2(\mathbb{C}^2,0)$ is the module of germs of holomorphic 2-forms on $(\mathbb{C}^2,0)$.

The $k$-th Milnor number of $\mathcal{F}$ and the $k$-th Tjurina number of
$\mathcal{F}$ along $C$ are then defined by
\[
\mu^k(\mathcal{F})
=\dim_{\mathbb{C}}
\frac{\mathcal{O}_2}
{\mathfrak{m}^k(P,Q)},
\qquad
\tau^k(\mathcal{F},C)
=\dim_{\mathbb{C}}
\frac{\mathcal{O}_2}
{\mathfrak{m}^k(P,Q)+(f)},
\]
respectively.
For $\omega=df$, we have $(P,Q)=J(f)$, so $\mu^k(\mathcal{F})=\mu^k(f)$
and $\tau^k(\mathcal{F},C)=\tau^k(f)$.

The following corollary relates the $k$-th Milnor number of a plane curve to its Milnor number.

\begin{corollary}[{\cite[Corollary 4.1]{Marcela}}]\label{mk}
Let $f\in\mathfrak{m}\subset\mathcal{O}_2$
 define an isolated singularity at
the origin $0\in\mathbb{C}^2$, that is, $\mu(f)<\infty$,
and let $m=\nu(f)$ be the order of $f$.
Then
\[\mu^k(f)=\begin{cases}
\mu(f)+k(k+1), & \text{if $0\leq k< m$,}\\
\mu(f)+k(k+1)-\frac{(k-m+2)(k-m+1)}{2}, & \text{if $k\geq m$.}
\end{cases}\]
\end{corollary}

Applying Corollary \ref{mk} to the product $fg$, we obtain

\begin{corollary}\label{coropro}
Let $f,g \in\mathfrak{m}\subset\mathcal{O}_2$ be reduced and coprime and let $m=\nu(f)$ and $\ell=\nu(g)$ be the orders of $f$ and $g$. Then
\[\mu^k(fg)=\begin{cases}
\mu(fg)+k(k+1), & \text{if $0\leq k< m+\ell$,}\\
\mu(fg)+k(k+1)-\frac{(k-(m+\ell)+2)(k-(m+\ell)+1)}{2}, & \text{if $k\geq m+\ell$.}
\end{cases}\]
\end{corollary}
\begin{proof}
Since $f$ and $g$ are reduced and coprime, $fg$ is reduced, so it has
an isolated singularity at the origin and $\mu(fg)<\infty$.
Since $\mathcal{O}_2$ is a domain, $\nu(fg)=\nu(f)+\nu(g)=m+\ell$.
Now apply Corollary~\ref{mk} to $fg$.
\end{proof}


\section{On a conjecture of Hussain--Yau--Zuo}\label{cociente-b}

In~\cite{Almiron}, P. Almir\'on studied the quotient between the Milnor number and
the Tjurina number of an isolated complete intersection singularity and
established its connections with other problems in singularity theory. As part
of this study, he proposed the following conjecture~\cite[Conjecture 4.10]{Almiron}.

\begin{conjecture}\label{conj:Almiron-b}
Let $f\in\mathfrak{m}\subset\mathcal{O}_3$ define an isolated hypersurface
singularity $(V,0)\subset(\C^3,0)$. Then
\[\frac{\mu(f)}{\tau(f)}<\frac{3}{2}.\]
\end{conjecture}

Almir\'on showed that Conjecture \ref{conj:Almiron-b} follows from Durfee's
conjecture $6 p_g\leq\mu$ (here $p_g$ denotes the \textit{geometric genus} of an isolated surface singularity), and verified it for several families of surface
singularities for which Durfee's conjecture is known to hold. Recently,
Ma--Zuo \cite{Ma-Zuo} proved Conjecture \ref{conj:Almiron-b} in full
generality, without relying on Durfee's conjecture, as a consequence of a sharp
upper bound for $\mu/\tau$ valid for isolated hypersurface singularities in any
dimension \cite[Theorem~4.5 and Corollary 4.8]{Ma-Zuo}.

In the context of the $k$-th Milnor and Tjurina numbers
of hypersurface singularities, Hussain--Yau--Zuo \cite[Conjecture~1.4]{Hussain-Yau}
conjectured the following inequality.

\begin{conjecture}\label{conj:k-b}
Let $f\in\mathfrak{m}\subset\mathcal{O}_3$ define an isolated hypersurface
singularity $(V,0)\subset(\C^3,0)$. Then for each integer $k\geq 0$,
\[\dfrac{\mu^k(f)}{\tau^k(f)}<\dfrac{3}{2}.\]
\end{conjecture}

Since $\mu^0(f)=\mu(f)$ and $\tau^0(f)=\tau(f)$, the case $k=0$ of
Conjecture~\ref{conj:k-b} is precisely Conjecture~\ref{conj:Almiron-b}, and hence
it holds by~\cite[Corollary~4.8]{Ma-Zuo}. The question is therefore whether the
bound $3/2$ persists for $k\geq 1$.
Hussain--Yau--Zuo considered the hypersurfaces defined by
$f_1=x^{3m+1}+y^{3m+1}+z^{3m+1}+x^{m+1}y^{m+1}z^{m+1}$ and
$f_2=x^{3m-1}+y^{3m-1}+z^{3m-1}+x^{m}y^{m}z^{m}$, with $m\geq 1$,
and verified explicitly that $\mu^{k}/\tau^{k}<3/2$ for $k=0,1,2$.
We show that Conjecture~\ref{conj:k-b} fails already for the simplest
isolated hypersurface singularity.

The computation of the $k$-th Tjurina numbers below rests on the following lemma,
which gives $\dim_{\mathbb{C}}\mathcal{O}_n/((f)+\mathfrak{m}^N)$.
 For a nonzero germ $h\in\mathcal{O}_2$, we denote by
$\mathrm{in}(h)$ its \emph{initial form}, that is, the nonzero
homogeneous component of lowest degree in its Taylor expansion.
Thus $\deg\mathrm{in}(h)=\nu(h)$.
\begin{lemma}\label{lem:quot}
Let $f\in\mathfrak{m}\subset\mathcal{O}_n$, let $\nu=\nu(f)$ and let $N\geq 0$ be an
integer. Put $\mathfrak{m}^j=\mathcal{O}_n$ for $j\leq 0$. Then
$\dim_{\mathbb{C}}\mathcal{O}_n/\mathfrak{m}^j=\binom{j+n-1}{n}$ for every $j\geq 0$,
$(\mathfrak{m}^N:f)=\mathfrak{m}^{N-\nu}$ and
\[
\dim_{\mathbb{C}}\frac{\mathcal{O}_n}{(f)+\mathfrak{m}^N}
=\dim_{\mathbb{C}}\frac{\mathcal{O}_n}{\mathfrak{m}^N}
-\dim_{\mathbb{C}}\frac{\mathcal{O}_n}{\mathfrak{m}^{N-\nu}}.
\]
In particular, for $N\geq\nu$ this equals
$\binom{N+n-1}{n}-\binom{N-\nu+n-1}{n}$.
\end{lemma}

\begin{proof}
For $j\geq 0$, the monomials of degree less than $j$ form a basis of
$\mathcal{O}_n/\mathfrak{m}^j$, and there are $\binom{d+n-1}{n-1}$ monomials of degree $d$;
hence $\dim\mathcal{O}_n/\mathfrak{m}^j=\sum_{d=0}^{j-1}\binom{d+n-1}{n-1}=\binom{j+n-1}{n}$.

Let $\varphi\colon\mathcal{O}_n\to\mathcal{O}_n/\mathfrak{m}^N$,
$\varphi(h)=fh+\mathfrak{m}^N$. It is $\mathbb{C}$-linear and its image is
$((f)+\mathfrak{m}^N)/\mathfrak{m}^N$. Its kernel is
$(\mathfrak{m}^N:f)=\{h : fh\in\mathfrak{m}^N\}$, so by the first isomorphism theorem
the image is isomorphic to $\mathcal{O}_n/(\mathfrak{m}^N:f)$.
The initial form of $fh$ is the product of the initial forms of $f$ and $h$, which is
nonzero; hence $\nu(fh)=\nu+\nu(h)$. Then $fh\in\mathfrak{m}^N$ if and only if
$\nu(h)\geq N-\nu$, that is, $(\mathfrak{m}^N:f)=\mathfrak{m}^{N-\nu}$.
Finally, let $\iota$ be the inclusion
$((f)+\mathfrak{m}^N)/\mathfrak{m}^N\hookrightarrow\mathcal{O}_n/\mathfrak{m}^N$ and let
$\pi$ be the canonical projection
$\mathcal{O}_n/\mathfrak{m}^N\to\mathcal{O}_n/((f)+\mathfrak{m}^N)$, which is well
defined since $\mathfrak{m}^N\subset(f)+\mathfrak{m}^N$. Then $\iota$ is injective,
$\pi$ is surjective and $\ker\pi=\operatorname{im}\iota$, so the sequence
$0\to((f)+\mathfrak{m}^N)/\mathfrak{m}^N\xrightarrow{\iota}\mathcal{O}_n/\mathfrak{m}^N
\xrightarrow{\pi}\mathcal{O}_n/((f)+\mathfrak{m}^N)\to 0$
is exact. It gives
$\dim\mathcal{O}_n/((f)+\mathfrak{m}^N)=\dim\mathcal{O}_n/\mathfrak{m}^N-
\dim\mathcal{O}_n/(\mathfrak{m}^N:f)$, and the stated formulas follow.
\end{proof}

\begin{proposition}\label{prop:A1}
Let $f=x^2+y^2+z^2\in\mathbb{C}\{x,y,z\}$. Then for every integer $k\geq 0$
\[
\mu^k(f)=\binom{k+3}{3},
\qquad
\tau^k(f)=(k+1)^2,
\]
and $\mu^k(f)/\tau^k(f)>3/2$ if and only if $k\geq 5$. In particular,
$\mu^5(f)=56$, $\tau^5(f)=36$ and $\mu^5(f)/\tau^5(f)=14/9>3/2$.
\end{proposition}

\begin{proof}
Since $J(f)=(2x,2y,2z)=\mathfrak{m}$, we have $\mathfrak{m}^kJ(f)=\mathfrak{m}^{k+1}$. Hence
$\mu^k(f)=\dim_{\mathbb{C}}\mathcal{O}_3/\mathfrak{m}^{k+1}=\binom{k+3}{3}$ by the formula
for $\dim_{\mathbb{C}}\mathcal{O}_n/\mathfrak{m}^j$ in Lemma~\ref{lem:quot} with $n=3$, $j=k+1$.

Put $N=k+1$. Since $\nu(f)=2$, Lemma~\ref{lem:quot} gives
\[
\tau^k(f)=\dim_{\mathbb{C}}\frac{\mathcal{O}_3}{\mathfrak{m}^{k+1}}
-\dim_{\mathbb{C}}\frac{\mathcal{O}_3}{\mathfrak{m}^{k-1}}
=\binom{k+3}{3}-\binom{k+1}{3}
\]
(for $k=0$ we have $\mathfrak{m}^{-1}=\mathcal{O}_3$ and $\binom{1}{3}=0$). Since
$(k+3)(k+2)(k+1)-(k+1)k(k-1)=6(k+1)^2$, we get $\tau^k(f)=(k+1)^2$.

Finally, $\mu^k(f)/\tau^k(f)>3/2$ is equivalent to
$2\mu^k(f)>3\tau^k(f)$, which can be rewritten as $(k+3)(k+2)(k+1)>9(k+1)^2$. Dividing by $k+1>0$,
this becomes $(k+3)(k+2)>9(k+1)$, that is, to $k^2-4k-3>0$.
The roots of this polynomial are $2\pm\sqrt{7}$ and its leading coefficient is positive,
so for $k\geq 0$ the inequality holds exactly when $k>2+\sqrt{7}$, that is, for the
integers $k\geq 5$, since $2<\sqrt{7}<3$.
\end{proof}

\begin{theorem}\label{thm:growth}
Let $n\geq 2$ and let $f\in\mathfrak{m}\subset\mathcal{O}_n=\mathbb{C}\{x_1,\ldots,x_n\}$
with $\mu(f)<\infty$. Then
$\displaystyle\lim_{k\to\infty}\mu^k(f)/\tau^k(f)=\infty$. In particular, for $n=3$
Conjecture~\ref{conj:k-b} fails for every such $f$ and every $k$ large enough.
\end{theorem}

\begin{proof}
Let $\nu=\nu(f)$ and $r=\mu(f)$. If $\nu(f)=1$, a biholomorphic change of coordinates reduces to
$f=x_1$. For every $k\geq1$, we then have
\[
\mu^k(f)=\binom{k+n-1}{n},\qquad
\tau^k(f)=\binom{k+n-2}{n-1},
\]
and hence
\[
\frac{\mu^k(f)}{\tau^k(f)}=\frac{k+n-1}{n}
\longrightarrow\infty.
\]
We may therefore assume that $\nu(f)\geq2$.
Write $f=\sum_{d\geq\nu}f_d$ with $f_d$
homogeneous of degree $d$. Then $\partial f/\partial x_i=\sum_{d\geq\nu}
\partial f_d/\partial x_i$, and each $\partial f_d/\partial x_i$ is zero or homogeneous
of degree $d-1\geq\nu-1$. Hence $\partial f/\partial x_i\in\mathfrak{m}^{\nu-1}$ for
every $i$, so $J(f)\subset\mathfrak{m}^{\nu-1}$ and
$\mathfrak{m}^kJ(f)\subset\mathfrak{m}^{k+\nu-1}$. Therefore
$\mu^k(f)=\dim_{\mathbb{C}}\mathcal{O}_n/\mathfrak{m}^kJ(f)\geq
\dim_{\mathbb{C}}\mathcal{O}_n/\mathfrak{m}^{k+\nu-1}$.
By the formula for $\dim_{\mathbb{C}}\mathcal{O}_n/\mathfrak{m}^j$ in
Lemma~\ref{lem:quot} with $j=k+\nu-1\geq 0$,
\[
\dim_{\mathbb{C}}\mathcal{O}_n/\mathfrak{m}^{k+\nu-1}
=\binom{(k+\nu-1)+n-1}{n}=\binom{k+\nu+n-2}{n}.
\]

Recall that $r=\mu(f)=\dim_{\mathbb{C}}\mathcal{O}_n/J(f)$. The algebra
$A=\mathcal{O}_n/J(f)$ is local, with maximal ideal $\overline{\mathfrak{m}}=\mathfrak{m}/J(f)$.
If $\overline{\mathfrak{m}}^i=\overline{\mathfrak{m}}^{i+1}$ for some $i$, then
$\overline{\mathfrak{m}}^i=0$ by Nakayama's lemma, since $\overline{\mathfrak{m}}^i$ is a
finitely generated $A$-module and $\overline{\mathfrak{m}}\cdot\overline{\mathfrak{m}}^i=
\overline{\mathfrak{m}}^i$. Hence, as long as $\overline{\mathfrak{m}}^i\neq 0$, the
inclusion $\overline{\mathfrak{m}}^{i+1}\subset\overline{\mathfrak{m}}^i$ is strict, so it
loses at least one dimension. Since $\dim_{\mathbb{C}}A=r$, this gives
$\dim_{\mathbb{C}}\overline{\mathfrak{m}}^i\leq r-i$ and, in particular,
$\overline{\mathfrak{m}}^r=0$, that is, $\mathfrak{m}^r\subset J(f)$. Therefore
$\mathfrak{m}^{k+r}\subset\mathfrak{m}^kJ(f)$, so $(f)+\mathfrak{m}^{k+r}\subset
(f)+\mathfrak{m}^kJ(f)$ and, for $N=k+r$,
$\tau^k(f)\leq\dim_{\mathbb{C}}\mathcal{O}_n/((f)+\mathfrak{m}^N)$.
Let $k\geq\nu$, so that $N\geq\nu$. By Lemma~\ref{lem:quot},
\begin{eqnarray*}
\dim_{\mathbb{C}}\frac{\mathcal{O}_n}{(f)+\mathfrak{m}^N}
=\binom{N+n-1}{n}-\binom{N-\nu+n-1}{n}
&=&\sum_{i=0}^{\nu-1}\binom{N+n-2-i}{n-1}\\
& &
\leq\nu\binom{N+n-2}{n-1},
\end{eqnarray*}
where the equality of sums follows from Pascal's rule
$\binom{a}{n}-\binom{a-1}{n}=\binom{a-1}{n-1}$ applied $\nu$ times.
Thus, for $k\geq\nu$,
\[
\frac{\mu^k(f)}{\tau^k(f)}\geq
\frac{\binom{k+\nu+n-2}{n}}{\nu\binom{k+r+n-2}{n-1}}.
\]
The numerator is a polynomial of degree $n$ in $k$ and the denominator a
polynomial of degree $n-1$, so the right-hand side tends to $\infty$.
\end{proof}

\begin{remark}
The plane curve conjecture was already disproved
in \cite[Section 7]{Marcela}.
The node $x^2+y^2$ provides the following simpler counterexample.
The same proof as
in Proposition~\ref{prop:A1}, now in two variables
(Lemma~\ref{lem:quot} with $n=2$), gives
$\mu^k(x^2+y^2)=\binom{k+2}{2}$ and $\tau^k(x^2+y^2)=2k+1$, so
$\mu^k/\tau^k>4/3$ if and only if $3k^2-7k-2>0$, that is, $k\geq 3$
(for instance $\mu^3=10$, $\tau^3=7$). Hence the analogous conjecture for
plane curves \cite[Conjecture~1.3]{Hussain-Yau} also fails.
\end{remark}

Hussain--Yau--Zuo supported Conjecture~\ref{conj:k-b} with the families
$f_1=x^{3m+1}+y^{3m+1}+z^{3m+1}+x^{m+1}y^{m+1}z^{m+1}$ and
$f_2=x^{3m-1}+y^{3m-1}+z^{3m-1}+x^{m}y^{m}z^{m}$, with $m\geq 1$:
they computed $\mu^k/\tau^k$ only for $k=0,1,2$, obtained values close to
$27/19<3/2$ for $m$ large, and wrote that they had no example with $k\geq 1$
and quotient close to $3/2$ \cite{Hussain-Yau}. Proposition~\ref{prop:A1} and
Theorem~\ref{thm:growth} show that checking only $k\leq 2$ cannot detect the failure of
Conjecture~\ref{conj:k-b}: for $x^2+y^2+z^2$ the quotient exceeds $3/2$ only from $k=5$ on,
and in general it tends to infinity with $k$.
Theorem~\ref{thm:growth} applies to $f_1$ and $f_2$, which have isolated singularities,
as required in Conjecture~\ref{conj:k-b}. Hence
$\mu^k(f_i)/\tau^k(f_i)\to\infty$ as $k\to\infty$ for every $m\geq 1$.
For instance, for $m=2$ and $f_1=x^7+y^7+z^7+x^3y^3z^3$, computations with
\textsc{Singular}~\cite{Singular} give $\mu^{16}(f_1)=2024$ and $\tau^{16}(f_1)=1344$,
so $\mu^{16}(f_1)/\tau^{16}(f_1)=253/168>3/2$.

\section{On the \texorpdfstring{$k$}{k}-th Milnor number of a
semiquasihomogeneous function}\label{k milnor}

A polynomial $f\in\mathbb{C}[x_1,\dots,x_n]$
is called \emph{quasihomogeneous} of degree
$d$ if there exists a vector
$w=(w_1,\dots,w_n)\in\mathbb{N}_{>0}^n$
and $d\in\mathbb{N}$ such that every
monomial $x^\alpha$ of $f$ with nonzero
coefficient satisfies
$\sum_{i=1}^{n}\alpha_i w_i=d$. In this
case we say that $f$ is of type
$(w_1,\dots,w_n,d)$. We call
$w_1,\dots,w_n$ the weights of the
variables and $d$ the weighted degree of
$f$, denoted as $\deg_w(f)=d$. A
quasihomogeneous polynomial $f$ is said to
be \textit{nondegenerate} if $0$ is an
isolated critical point. A power series
$f\in\mathcal{O}_n$ is said to be
\textit{semiquasihomogeneous} (SQH) of
degree $d$ with weights $w_1,\ldots,w_n$ if
it is of the form $f=f_w+f'$, where $f_w$
is a nondegenerate quasihomogeneous
polynomial of degree $d$ with weights $w_i$
and $f'$ is a power series all of whose
monomials have weighted degree strictly
greater than $d$.

Let $K[[x_1,\ldots,x_n]]$ be the ring of
formal power series over a field $K$. For
brevity, we write $K[[x_1,\ldots,x_n]]$ as
$K[[x]]$. Let $\mathfrak{m}$ be the maximal
ideal of $K[[x]]$. In this situation, we
can also define the $k$-th Milnor number of
$f\in K[[x]]$ as
$\mu^k(f):=\dim_K K[[x]]/\mathfrak{m}^kJ(f)$,
where $J(f)$ is the Jacobian ideal of $f$ in
$K[[x]]$. If the dimension is infinite, we
simply write $\mu^k(f)=\infty$. In what
follows in this section $K=\mathbb{C}$.

For a SQH power series $f$ of weighted
degree $d$ with weights $w_1,\dots,w_n$ we
have $\mu(f)=\mu(f_w)<\infty$. This is proved
in \cite[p.~194]{Arnold} for convergent power
series, by a topological argument (the degree
of the gradient map) which does not apply to
formal series. For formal power series,
\cite[Proposition~2]{BGM} gives, by an
algebraic argument,
$\mu(f)=\prod_{i=1}^{n}(d/w_i-1)$, which
depends only on $d$ and $w$; applying it to
$f$ and to $f_w$ gives $\mu(f)=\mu(f_w)$.
Hence $\mathfrak{m}^N\subset J(f)$ for some
$N$ and $\mathfrak{m}^kJ(f)\supset
\mathfrak{m}^{k+N}$. Since
$\mathbb{C}\{x\}/\mathfrak{m}^s=
\mathbb{C}[[x]]/\mathfrak{m}^s$ for every
$s\geq 0$, the number $\mu^k(f)$ is the same
whether computed in $\mathbb{C}\{x\}$ or in
$\mathbb{C}[[x]]$.

In \cite[Conjecture~4.65]{Shi},
Shi--Yau--Zuo proposed the following
conjecture.

\begin{conjecture}\label{conjecture Shi}
Let $f=f_w+f'\in\mathbb{C}[[x]]$ be a
semiquasihomogeneous $\operatorname{(SQH)}$
power series with respect to a weight vector
$w\in\mathbb{N}_{>0}^n$, and assume that
$\mu(f_w)<\infty$. Then
$\mu^k(f_w)=\mu^k(f)$ for every integer
$k\geq 0$.
\end{conjecture}

In~\cite[Theorem~C]{Shi}, the authors proved
Conjecture~\ref{conjecture Shi} under the
additional assumption that
$k\leq\min_i\{\nu(f_{x_i})\}$.

Liu--Yau--Zhang
\cite[Main Theorem~A]{Liu-Yau-Zhang}
recently proved
Conjecture~\ref{conjecture Shi} for
convergent germs in the equal-weight case,
for every $k\geq 0$ and every $n\geq 2$.
More precisely, their result applies to
germs $f=f_r+f_{r+1}+\cdots$ with
$r\geq 2$, where each $f_j$ is homogeneous
of degree $j$ and $V_+(f_r)\subset
\mathbb{P}^{n-1}$ is smooth (equivalently,
$f_r$ has an isolated critical point at the
origin).

In this section, we show that the conjecture
fails in dimension three for unequal weights.
Our counterexample lies outside the range
$k\leq\min_i\{\nu(f_{x_i})\}$ established in
\cite[Theorem~C]{Shi}. We also prove that the
conjecture holds in dimension two over
$\mathbb{C}$ for arbitrary positive weights
and every $k\geq 0$; see
Theorem~\ref{Teorema-semi}.

To this end, consider the germ
\[
f=x^3+y^7+z^7+xy^2z^3
\]
in $\mathcal{O}_3=\mathbb{C}\{x,y,z\}$. Its
principal part with respect to the weight
vector $w=(7,3,3)$ is
\[
f_w=x^3+y^7+z^7.
\]
Since $f_w$ has weighted degree $21$, the
monomial $xy^2z^3$ has weighted degree $22$,
and $\mu(f_w)=72<\infty$, the germ $f$ is
semiquasihomogeneous.
Using \textsc{Singular}~\cite{Singular}
with a local monomial ordering, we compute
the dimensions of $\mathcal{O}_3/J(h)$ and
$\mathcal{O}_3/\mathfrak{m}^6J(h)$
for $h=f_w$ and $h=f$, and we obtain
$\mu(f_w)=\mu(f)=72$ and
$\mu^6(f_w)=200\neq 199=\mu^6(f)$,
which gives a negative answer to
Conjecture~\ref{conjecture Shi}
in dimension three.

Note that $\nu(f_{x})=2$ and
$\nu(f_{y})=\nu(f_{z})=5$, so
$k=6>\min\{\nu(f_{x}),\nu(f_{y}),
\nu(f_{z})\}=2$, which is consistent with
\cite[Theorem~C]{Shi}.

We now prove that the conjecture holds in
dimension two over the ring
$\mathbb{C}[[x,y]]$.

In the proof of the next theorem we use
\cite[Corollary~4.1]{Marcela}, which gives
the formula for $\mu^k(f)$ of a plane
curve. It is stated for convergent power
series, but its proof uses only algebraic
arguments (orders of products, relatively
prime elements of a unique factorization
domain, and exact sequences of vector
spaces), which are also valid in
$\mathbb{C}[[x,y]]$. Hence the formula
holds for every $f\in\mathbb{C}[[x,y]]$
with $\mu(f)<\infty$.

\begin{theorem}\label{Teorema-semi}
Suppose that $f=f_w+f'\in\mathbb{C}[[x,y]]$
is a semiquasihomogeneous
$\operatorname{(SQH)}$ power series with
respect to a weight vector
$w\in\mathbb{N}_{>0}^2$. Then
$\mu^k(f_w)=\mu^k(f)$
for every integer $k\geq 0$.
\end{theorem}
\begin{proof}
By \cite[Proposition~2]{BGM}, applied to $f$
and to $f_w$, which are SQH with respect to
the same $w$ and the same degree $d$,
$\mu(f)=\prod_{i=1}^{2}(d/w_i-1)=\mu(f_w)
<\infty$. Hence $f$ and $f_w$ have isolated
singularities, and by
Corollary~\ref{mk} the number $\mu^k(g)$,
for $g=f$ or $g=f_w$, depends only on $k$,
$\mu(g)$ and $\nu(g)$.

We claim that $\nu(f_w)=\nu(f)$. After
exchanging the variables if necessary,
assume that $w_1\geq w_2$. Since $f_w$ has
an isolated critical point, its gradient
cannot vanish identically on the $x$-axis.
Hence $f_w$ contains a monomial $x^Ay^c$
with $c\in\{0,1\}$. Its weighted degree
satisfies
$d=Aw_1+cw_2>(A+c-1)w_1$.
Every monomial of total degree at most
$A+c-1$ has weighted degree at most
$(A+c-1)w_1<d$, so no such monomial occurs
in either $f_w$ or $f'$. Moreover, the
monomial $x^Ay^c$ cannot be cancelled by
$f'$, whose monomials have weighted degree
strictly greater than $d$. Thus
$\nu(f_w)=\nu(f)=A+c$.

Since $\mu(f_w)=\mu(f)$ and
$\nu(f_w)=\nu(f)$, we conclude that
$\mu^k(f_w)=\mu^k(f)$ for every integer
$k\geq 0$.
\end{proof}

The following example illustrates
Theorem~\ref{Teorema-semi}.
\begin{example}
Consider the function
$f(x,y)=x^3+y^4+x^2y^3$. Its
semiquasihomogeneous principal part with
respect to the weight vector $w=(4,3)$ is
$f_w(x,y)=x^3+y^4$ and $x^2y^3$ has
weighted degree $17>12$. Here
$\nu(f)=\nu(f_w)=3$ and, since
$J(f_w)=(x^2,y^3)$,
$\mu(f)=\mu(f_w)=
\dim_{\mathbb{C}}\mathcal{O}_2/(x^2,y^3)=6$.
Using Corollary~\ref{mk}, we get
\[
\mu^k(f)=\mu^k(f_w)=
\begin{cases}
6, & k=0,\\[2mm]
\dfrac{k^2+5k+10}{2}, & k\geq1.
\end{cases}
\]
Indeed, for $k=0$ we get $\mu^0=6$. For
$k=1,2$ (that is, $k<3$)
Corollary~\ref{mk} gives
$6+k(k+1)=(k^2+5k+10)/2$, and for
$k\geq 3$ it gives
$6+k(k+1)-(k-1)(k-2)/2=(k^2+5k+10)/2$.
\end{example}

\section{The \texorpdfstring{$k$}{k}-th Milnor and Tjurina numbers of a pair of
holomorphic function germs}\label{k pares}

Let $f,g\in \mathfrak{m}\subset\mathcal{O}_2$ and put $P=gf_x-fg_x$ and $Q=gf_y-fg_y$.
The Milnor number of the pair $(f,g)$ is
$\mu(f,g)=\dim_{\mathbb{C}}\mathcal{O}_2/(P,Q)\in\mathbb{N}\cup\{\infty\}$.
When $\mu(f,g)<\infty$, we denote by $\mathcal F_{f/g}$ the
foliation defined by $\omega_{f/g}=g\,df-f\,dg$.

For each integer $k\geq 0$, we define the $k$-th Milnor number of the pair $(f,g)$ by
\[
\mu^{k}(f,g)=\dim_{\mathbb C} \dfrac{\mathcal{O}_2}{\mathfrak{m}^k(gf_x-fg_x,gf_y-fg_y)}\in\mathbb{N}\cup\{\infty\}.
\]
Note that $\mu^0(f,g)=\mu(f,g)$.

For $f,g\in\mathcal{O}_2$, a point $[\alpha:\beta]\in\mathbb{P}^1_{\mathbb{C}}$ and an integer $k\geq 0$, the $k$-th Tjurina number of the pair $(f,g)$ with respect to $C_{[\alpha:\beta]}:\alpha f+\beta g=0$ is defined by 
\[
\tau^{k}_{[\alpha:\beta]}(f,g)=\dim_{\mathbb C}  \frac{\mathcal{O}_2}{\mathfrak{m}^k(gf_x-fg_x,gf_y-fg_y)+ (\alpha f+\beta g)}.
\]
For integers $k\geq 0$ and $m_0\geq 0$, put
$\varepsilon_k(m_0)=0$ if $0\leq k<m_0$ and
$\varepsilon_k(m_0)=\dfrac{(k-m_0+1)(k-m_0)}{2}$ if $k\geq m_0$.

\begin{theorem}\label{teo1}
Let $f,g\in\mathfrak{m}\subset\mathcal{O}_2$ be relatively prime reduced germs such that $\mu(f,g)<\infty$. Let $m_0=\min\{\nu(gf_x-fg_x),\nu(gf_y-fg_y)\}$. Then, for every integer $k\geq 0$,

   \[\mu^k(f,g)=\mu(f,g)+k(k+1)-\varepsilon_k(m_0),\]
that is,
   \[\mu^k(f,g)=\begin{cases}
			\mu(f,g)+k(k+1), & \text{if $0\leq k< m_0$,}\\
            \mu(f,g)+k(k+1)-\dfrac{(k-m_0+1)(k-m_0)}{2}, & \text{if $k\geq m_0$.}
		 \end{cases}\]
\end{theorem}
\begin{proof}
Since $\mu(f,g)<\infty$, the germs $P$ and $Q$ are relatively prime. Interchanging the names of $P$ and $Q$ if necessary (this does not change the ideal $(P,Q)$), we may assume that $\nu(P)\leq\nu(Q)$, so that $\nu(P)=m_0$.

Applying the formula from \cite[Lemma 3.3]{Marcela} to the relatively prime germs $P$ and $Q$, with $\nu(P)\leq\nu(Q)$, and $J=\mathfrak{m}^k$, we obtain

\begin{equation}\label{ecua}
\mu^{k}(\mathcal{F}_{f/g})
=\mu(\mathcal{F}_{f/g})
+\dim\frac{\mathcal{O}_{2}}{\mathfrak{m}^k}
+\dim\frac{\mathcal{O}_{2}}{\mathfrak{m}^{k}+(P)}.
\end{equation}

By definition, $\mu^k(\mathcal{F}_{f/g})=\mu^k(f,g)$ and $\mu(\mathcal{F}_{f/g})=\mu(f,g)$. Moreover, $\dim\mathcal{O}_2/\mathfrak{m}^k=k(k+1)/2$, the number of monomials of degree less than $k$.
Next, $(\mathfrak{m}^k+(P))/\mathfrak{m}^k\cong\mathcal{O}_2/(\mathfrak{m}^k:P)$. Since $\nu(Pa)=m_0+\nu(a)$ for every $a\in\mathcal{O}_2$, we have $(\mathfrak{m}^k:P)=\mathfrak{m}^{k-m_0}$ if $k\geq m_0$ and $(\mathfrak{m}^k:P)=\mathcal{O}_2$ if $k< m_0$. Hence
$\dim(\mathfrak{m}^k+(P))/\mathfrak{m}^k=\varepsilon_k(m_0)$ and
\[
\dim\frac{\mathcal{O}_{2}}{\mathfrak{m}^{k}+(P)}
=\frac{k(k+1)}{2}-\varepsilon_k(m_0)
=
\begin{cases}
\frac{k(k+1)}{2}, 
& \text{if $0\leq k< m_0$,}\\
\frac{k(k+1)}{2}
-\frac{(k-m_0+1)(k-m_0)}{2}, 
& \text{if $k\geq m_0$.}
\end{cases}
\]

Substituting these expressions into \eqref{ecua} completes the proof of the theorem.
\end{proof}

Theorem \ref{teo1} shows that $\mu^k(f,g)$ depends only on
$\mu(f,g)$ and $m_0$. Moreover, for $f,g\in\mathcal{O}_2$ and every integer $k\geq 0$, $\mu^k(f,g)$ is invariant
under biholomorphic changes of the source coordinates and
invertible linear changes of the generators. Indeed, if
\[
\begin{pmatrix}f_1\\ g_1\end{pmatrix}
=A\begin{pmatrix}f\circ\Phi\\ g\circ\Phi\end{pmatrix},
\qquad A\in\operatorname{GL}_2(\mathbb{C}),
\]
where $\Phi$ is a biholomorphic germ, then
\[
g_1\,df_1-f_1\,dg_1
=\det(A)\,\Phi^*(g\,df-f\,dg).
\]
The coefficients $P_1,Q_1$ of $g_1\,df_1-f_1\,dg_1$ are $\det(A)$ times the coefficients of $\Phi^*(P\,dx+Q\,dy)$, which are obtained from $P\circ\Phi$ and $Q\circ\Phi$ by multiplying by the Jacobian matrix of $\Phi$, invertible over $\mathcal{O}_2$. Hence $(P_1,Q_1)=(P\circ\Phi,Q\circ\Phi)=\Phi^*(P,Q)$.
As $\Phi^*$ is a ring automorphism of $\mathcal{O}_2$ with $\Phi^*(\mathfrak m)=\mathfrak m$,
we obtain $\mu^k(f_1,g_1)=\mu^k(f,g)$ for every integer $k\geq0$.

Using the formula for the $k$-th Tjurina numbers established in \cite[Theorem~8.2]{Marcela}, we relate $\tau^k_{[\alpha:\beta]}(f,g)$
to the $k$-th Tjurina number of the curve
$C_{[\alpha:\beta]}$.
 
\begin{proposition}\label{teo_2}
Let $f,g\in\mathfrak{m}\subset\mathcal{O}_2$ be relatively prime
reduced germs such that $\mu(f,g)<\infty$.
Let $[\alpha:\beta]\in\mathbb{P}^1_{\mathbb{C}}$
and let
$C_{[\alpha:\beta]}:\alpha f+\beta g=0$
be the corresponding member of the pencil.
Then $C_{[\alpha:\beta]}$ is reduced and, for
every integer $k\geq 0$,
$\tau^{k}_{[\alpha:\beta]}(f,g)=
\tau^k(C_{[\alpha:\beta]})+I_0(f,g)$.
\end{proposition}
\begin{proof}
Since $\mu(f,g)<\infty$, we have that
$\mathcal{F}_{f/g}$ is a foliation with an
isolated singularity at the origin
$0\in\mathbb{C}^2$. Put
$P=gf_x-fg_x$ and $Q=gf_y-fg_y$. Then
$\tau^k_{[\alpha:\beta]}(f,g)=
\tau^k(\mathcal{F}_{f/g},C_{[\alpha:\beta]})$
is the definition of the $k$-th Tjurina number
of $\mathcal{F}_{f/g}$ along
$C_{[\alpha:\beta]}$. Moreover,
$C_{[\alpha:\beta]}$ is
$\mathcal{F}_{f/g}$-invariant, since
$\omega_{f/g}\wedge d(\alpha f+\beta g)=(\alpha f+\beta g)\,df\wedge dg$ is
a multiple of $\alpha f+\beta g$.
Since $\mu(f,g)<\infty$, the curve
$C_{[\alpha:\beta]}$ is reduced by
\cite[Lemma~4.1]{FFMR}, which proves the
first assertion.
Then, it follows from
\cite[Theorem~8.2]{Marcela} applied to
$\mathcal{F}_{f/g}$ and
$C_{[\alpha:\beta]}:\alpha f+\beta g=0$ that
\begin{equation*}
\tau^k_{[\alpha:\beta]}(f,g)=
\tau^k(\mathcal{F}_{f/g},C_{[\alpha:\beta]})=
\tau^k(C_{[\alpha:\beta]})+
\gsv(\mathcal{F}_{f/g},C_{[\alpha:\beta]}),
\end{equation*}
where $\gsv(\mathcal{F}_{f/g},C_{[\alpha:\beta]})$
denotes the \gsv-index of $\mathcal{F}_{f/g}$
with respect to $C_{[\alpha:\beta]}$
(see \cite{GSV}). The proof ends by the equality
$\gsv(\mathcal{F}_{f/g},C_{[\alpha:\beta]})=
I_0(f,g)$, see \cite[Proposition 4.8]{FFMR}.
\end{proof}

In
particular, under the hypotheses of Proposition \ref{teo_2}, the difference
$\tau^{k}_{[\alpha:\beta]}(f,g)-
\tau^k(C_{[\alpha:\beta]})=I_0(f,g)$
is independent of both $k$ and the chosen reduced member of the pencil.

For $f,g\in\mathcal{O}_2$ such that $\mu(f,g)<\infty$, in \cite[(2.1) Lemma]{Suwa}, Suwa studied  the unfoldings of germs of meromorphic functions and introduced the invariant
\begin{equation}\label{def:eta}
\eta(f,g)=\dim_{\mathbb{C}}\frac{(f,g)}{(f_xg-g_xf,f_yg-g_yf)}.
\end{equation}

\begin{lemma}\label{lema-eta}
Let $f,g\in\mathcal{O}_2$ be such that $\mu(f,g)<\infty$. Then $\eta(f,g)$ and $I_0(f,g)$ are finite and
\begin{equation}\label{I0-mu-eta}
\mu(f,g)=\eta(f,g)+I_0(f,g).
\end{equation}
In particular, $I_0(f,g)=\mu(f,g)-\eta(f,g)$.
\end{lemma}
\begin{proof}
Put $P=gf_x-fg_x$ and $Q=gf_y-fg_y$.
Since $(P,Q)\subseteq(f,g)$, there is an exact sequence of vector spaces
\[
0\longrightarrow\frac{(f,g)}{(P,Q)}
\longrightarrow\frac{\mathcal O_2}{(P,Q)}
\longrightarrow\frac{\mathcal O_2}{(f,g)}
\longrightarrow0.
\]
Since $\mu(f,g)=\dim\mathcal{O}_2/(P,Q)<\infty$, the other two spaces have finite dimension, which are $\eta(f,g)$ and $I_0(f,g)$ respectively, and
$\mu(f,g)=\eta(f,g)+I_0(f,g)$.
\end{proof}

The following theorem relates $\mu^k(f,g)$ and
$\tau^k_{[\alpha:\beta]}(f,g)$.

\begin{theorem}\label{teo3}
Let $f,g\in\mathfrak{m}\subset\mathcal{O}_2$ be relatively
prime reduced germs such
that $\mu(f,g)<\infty$. Let $[\alpha:\beta]\in\mathbb{P}^1_{\mathbb{C}}$
and let
$C_{[\alpha:\beta]}:\alpha f+\beta g=0$
be the corresponding member of the pencil,
and let
$m_0=\min\{\nu(gf_x-fg_x),\nu(gf_y-fg_y)\}$.
Then, for every integer $k\geq 0$,
\[
\mu^k(f,g)-\tau^k_{[\alpha:\beta]}(f,g)=
\eta(f,g)+k(k+1)-\tau^{k}(C_{[\alpha:\beta]})
-\varepsilon_k(m_0).
\]
\end{theorem}
\begin{proof}
By Proposition \ref{teo_2} and Lemma \ref{lema-eta},
\[
\tau^k_{[\alpha:\beta]}(f,g)=
\tau^k(C_{[\alpha:\beta]})+I_0(f,g)=
\tau^k(C_{[\alpha:\beta]})+\mu(f,g)-\eta(f,g).
\]
By Theorem \ref{teo1},
$\mu^k(f,g)=\mu(f,g)+k(k+1)-\varepsilon_k(m_0)$.
Subtracting, $\mu(f,g)$ cancels and we obtain
the formula.
\end{proof}

The following example illustrates Theorem \ref{teo3}.

\begin{example}
Let $f(x,y)=x^2+y^3$ and $g(x,y)=xy$ in
$\mathcal{O}_{2}$. Consider the fiber
\[
C_{[\alpha:\beta]}:\alpha(x^2+y^3)+\beta xy=0,
\qquad [\alpha:\beta]\in\mathbb P^1_{\mathbb{C}}.
\]
The curve $C_{[1:0]}$ is the cusp, and for
$[\alpha:\beta]\neq[1:0]$ the curve
$C_{[\alpha:\beta]}$ is a node. In particular,
all of them are reduced. We have
$\tau^k(C_{[1:0]})=2k+2$, and for
$[\alpha:\beta]\neq[1:0]$,
$\tau^k(C_{[\alpha:\beta]})=2k+1$.
Moreover, $P=x^2y-y^4$ and $Q=2xy^3-x^3$
both have order $3$, and
\[
m_0=3,\quad \eta(f,g)=7,\quad
\mu(f,g)=12,\quad I_0(f,g)=5.
\]
By Theorem~\ref{teo1} (and checked with
\textsc{Singular}),
\[
\mu^k(f,g)=\begin{cases}
k^{2}+k+12, & \text{if } 0\leq k< 3,\\[1ex]
\dfrac{k^{2}+7k+18}{2}, & \text{if } k\geq 3.
\end{cases}
\]
By Proposition \ref{teo_2},
$\tau^k_{[\alpha:\beta]}(f,g)=2k+6$ for
$[\alpha:\beta]\neq[1:0]$, and
$\tau^k_{[1:0]}(f,g)=2k+7$. Therefore
\[
\mu^k(f,g)-\tau^k_{[\alpha:\beta]}(f,g)=
\begin{cases}
k^2-k+6, & 0\leq k< 3,\\[1ex]
\dfrac{k^{2}+3k+6}{2}, & k\geq 3,
\end{cases}\]
and 
\[\mu^k(f,g)-\tau^k_{[1:0]}(f,g)=
\begin{cases}
k^{2}-k+5, & 0\leq k< 3,\\[1ex]
\dfrac{k^{2}+3k+4}{2}, & k\geq 3.
\end{cases}
\]
respectively.
\end{example}

\section{Bifurcation formula for the pair}\label{bifu}

In this section we use
Theorem~\ref{teo1} and
Corollary~\ref{coropro} to compare the
$k$-th Milnor number $\mu^k(f,g)$ of the
pair with the $k$-th Milnor number
$\mu^k(fg)$ of the product, in the spirit
of the bifurcation formula of
\cite{FFMR}. We first compute the algebraic
multiplicity of the foliation $\mathcal{F}_{f/g}$
in terms of the bifurcation fibers of the pencil.

\noindent Let $f,g\in\mathfrak{m}\subset\mathcal{O}_2$ be relatively
prime germs such that $\mu(f,g)<\infty$.
Let $C_{[\alpha:\beta]}=\{\alpha f+\beta g=0\}$ denote
the member of the pencil $\mathcal{P}_{f,g}$ defined by
$[\alpha:\beta]\in\mathbb P^1_{\mathbb C}$.
By \cite[Lemma 4.1]{FFMR}, every $C_{[\alpha:\beta]}$ is reduced,
so $\mu(C_{[\alpha:\beta]})<\infty$. The
generic Milnor number and the generic multiplicity of the pencil are
\[
\mu_{\mathrm{gen}}
=\min_{[\alpha:\beta]\in\mathbb P^1_{\mathbb C}}
  \mu(C_{[\alpha:\beta]})
\quad\text{and}\quad
\nu_{\mathrm{gen}}
=\min_{[\alpha:\beta]\in\mathbb P^1_{\mathbb C}}
  \nu(\alpha f+\beta g).
\]
The bifurcation set is the finite set
(its finiteness is proved in \cite{FFMR})
\[
\mathcal B(f,g)=
\bigl\{[\alpha:\beta]\in\mathbb P^1_{\mathbb C}:
\mu(C_{[\alpha:\beta]})>\mu_{\mathrm{gen}}\bigr\},
\]
and we put
\[
\mathcal B^*(f,g)
=\mathcal B(f,g)\setminus\{[1:0],[0:1]\}.
\]
Note that $[1:0]$ and $[0:1]$ correspond to $f=0$ and $g=0$.
For $[\alpha:\beta]\in\mathcal B^*(f,g)$ we call
$C_{[\alpha:\beta]}$ a \emph{bifurcation fiber} of $\mathcal{P}_{f,g}$
other than $f$ and $g$.

\par The following proposition was established in \cite[Corollary 5.3]{Tjurina-par}; here we provide a different proof that does not use the bifurcation formula given in \cite{FFMR}.
\begin{proposition}\label{multfol}
Let $f,g\in\mathfrak{m}\subset\mathcal{O}_2$ be relatively
prime germs such that $\mu(f,g)<\infty$, and let
$h_1,\ldots,h_r$ be the distinct bifurcation fibers of
$\mathcal{P}_{f,g}$ other than $f$ and $g$. Then
\[
\nu(\mathcal{F}_{f/g})=\nu(f)+\nu(g)-1+\sum_{j=1}^r\bigl(\nu(h_j)-\nu_{\mathrm{gen}}\bigr).
\]
Moreover, at most one summand is nonzero, and this happens
if and only if $\nu(f)=\nu(g)$ and the initial forms of $f$ and $g$
are proportional.
\end{proposition}

\begin{proof}
Put $M=\nu(f)+\nu(g)$, $P=gf_x-fg_x$ and $Q=gf_y-fg_y$.
As in the proof of Theorem~\ref{kbifu} below,
$\nu(P),\nu(Q)\geq M-1$. Let
$P_{M-1}$ and $Q_{M-1}$ be the homogeneous parts of degree $M-1$ of
$P$ and $Q$. By Euler's identity,
\[
xP_{M-1}+yQ_{M-1}=\bigl(\nu(f)-\nu(g)\bigr)\,\mathrm{in}(f)\,\mathrm{in}(g).
\]

If $\nu(f)\neq\nu(g)$, then $(P_{M-1},Q_{M-1})\neq(0,0)$, so
$\nu(\mathcal{F}_{f/g})=M-1$. If $\nu(f)=\nu(g)$ and the initial forms
of $f$ and $g$ are not proportional, then $P_{M-1}=Q_{M-1}=0$ would give
$\mathrm{in}(g)\,d\,\mathrm{in}(f)-\mathrm{in}(f)\,d\,\mathrm{in}(g)=0$, that is,
$d\bigl(\mathrm{in}(f)/\mathrm{in}(g)\bigr)=0$, a contradiction; hence
again $\nu(\mathcal{F}_{f/g})=M-1$. In both cases
$\mathrm{in}(\alpha f+\beta g)\neq 0$ has degree $\min\{\nu(f),\nu(g)\}$
for every $[\alpha:\beta]\notin\{[1:0],[0:1]\}$, so every member of
$\mathcal{P}_{f,g}$ other than $f$ and $g$ has multiplicity
$\nu_{\mathrm{gen}}$, and the sum vanishes.

Assume now that $\nu(f)=\nu(g)$ and $\mathrm{in}(f)=c\,\mathrm{in}(g)$
with $c\in\mathbb{C}^*$. Then $\nu_{\mathrm{gen}}=\nu(f)=\nu(g)$, and
$h=f-cg$ satisfies $\nu(h)>\nu_{\mathrm{gen}}$. It is the only member of
$\mathcal{P}_{f,g}$ with this property: two distinct such members would
generate $\mathcal{P}_{f,g}$, and then every member would have
multiplicity greater than $\nu_{\mathrm{gen}}$. Since
$g\,dh-h\,dg=g\,df-f\,dg$ and $\nu(h)\neq\nu(g)$,
the case $\nu(f)\neq\nu(g)$ applied
to the pair $(h,g)$ gives
\[
\nu(\mathcal{F}_{f/g})=\nu(h)+\nu(g)-1
=\nu(f)+\nu(g)-1+\bigl(\nu(h)-\nu_{\mathrm{gen}}\bigr).
\]
It remains to show that $h$ is a bifurcation fiber. Suppose that
$\mu(h)=\mu_{\mathrm{gen}}$. Since $\mathcal B(f,g)$ is finite, the family
$h_s=f-(c+s)g$, $s\in(\mathbb{C},0)$, consists of reduced curves with
$\mu(h_s)=\mu_{\mathrm{gen}}$ for every $s$. By the L\^e--Ramanujam
theorem \cite{LR}, this $\mu$-constant family is topologically trivial,
and since the multiplicity of a reduced plane curve germ is a topological
invariant (see \cite{GLS}), all its members have the same multiplicity. This contradicts
$\nu(h_s)=\nu_{\mathrm{gen}}<\nu(h)=\nu(h_0)$ for $s\neq 0$. Hence
$\mu(h)>\mu_{\mathrm{gen}}$, so $h\in\{h_1,\ldots,h_r\}$, while the other
$h_j$ have multiplicity $\nu_{\mathrm{gen}}$.
\end{proof}

\begin{theorem}\label{kbifu}
Let $f,g\in\mathfrak{m}\subset\mathcal{O}_2$ be relatively
prime reduced germs
such that $\mu(f,g)<\infty$. Let
$m=\nu(f)$, $\ell=\nu(g)$, $M=m+\ell$ and
$m_0=\min\{\nu(gf_x-fg_x),\nu(gf_y-fg_y)\}$,
and put
\[
\Delta=\sum_{[\alpha:\beta]\in
\mathcal{B}^*(f,g)}
\bigl(\mu(C_{[\alpha:\beta]})
-\mu_{\mathrm{gen}}\bigr)
\quad\text{and}\quad
\delta=\sum_{[\alpha:\beta]\in
\mathcal{B}^*(f,g)}
\bigl(\nu(\alpha f+\beta g)
-\nu_{\mathrm{gen}}\bigr).
\]
Then $m_0=M-1+\delta\geq M-1$ and, for every integer $k\geq 0$,
\[
\mu^k(f,g)=\mu^k(fg)+\Delta+\varepsilon_k(M-1)-\varepsilon_k(M-1+\delta).
\]
\end{theorem}

\begin{proof}
Put $P=gf_x-fg_x$ and $Q=gf_y-fg_y$.
Since $\nu(f_x),\nu(f_y)\geq m-1$ and $\nu(g_x),\nu(g_y)\geq \ell-1$,
we have $\nu(gf_x),\nu(fg_x),\nu(gf_y),\nu(fg_y)\geq M-1$. Hence
$\nu(P)\geq M-1$ and $\nu(Q)\geq M-1$, that is, $m_0\geq M-1$.
Since $m_0=\nu(\mathcal{F}_{f/g})$, Proposition~\ref{multfol} gives
$m_0=M-1+\delta$.

By the bifurcation formula of
\cite[Theorem 4.6]{FFMR},
$\mu(f,g)=\mu(fg)+\Delta$. Since $f$ and $g$ are reduced and relatively prime, Corollary~\ref{coropro} gives
$\mu^k(fg)=\mu(fg)+k(k+1)-\varepsilon_k(M-1)$.
By Theorem~\ref{teo1},
$\mu^k(f,g)=\mu(f,g)+k(k+1)-\varepsilon_k(m_0)$.
Subtracting, and using $\mu(f,g)-\mu(fg)=\Delta$,
\[
\mu^k(f,g)-\mu^k(fg)=\Delta+\varepsilon_k(M-1)-\varepsilon_k(m_0).
\]
\end{proof}

\begin{remark}\label{casesdelta}
By Proposition~\ref{multfol}, $\delta>0$ if and only if $m=\ell$ and
the initial forms of $f$ and $g$ are proportional; in this case
$\delta=\nu(h)-\nu_{\mathrm{gen}}$, where $h$ is the unique member of
$\mathcal{P}_{f,g}$ of multiplicity greater than $\nu_{\mathrm{gen}}$.
Substituting the definition of $\varepsilon_k$ in the formula of
Theorem~\ref{kbifu} gives, according to the value of $\delta$:
\begin{enumerate}
\item If $\delta=0$, that is, $m_0=M-1$, then
$\varepsilon_k(M-1)-\varepsilon_k(m_0)=0$, and so
\[
\mu^k(f,g)=\mu^k(fg)+\Delta
\quad\text{for every } k\geq 0.
\]
\item If $\delta=1$, that is, $m_0=M$, then $\varepsilon_k(M-1)=\varepsilon_k(M)=0$ for $k<M$, and
$\varepsilon_k(M-1)-\varepsilon_k(M)=k-M+1$ for $k\geq M$. Hence
\[
\mu^k(f,g)=
\begin{cases}
\mu^k(fg)+\Delta,
& \text{if } 0\leq k<M,\\[1ex]
\mu^k(fg)+\Delta+(k-M+1),
& \text{if } k\geq M.
\end{cases}
\]
\item If $\delta\geq 2$, that is, $m_0>M$, then $\varepsilon_k(M-1)=\varepsilon_k(m_0)=0$ for $k<M$;
for $M\leq k<m_0$, $\varepsilon_k(m_0)=0$ and
$\varepsilon_k(M-1)=\frac{(k-M+2)(k-M+1)}{2}$; for $k\geq m_0$ both are given
by their second expression. Hence
\[
\mu^k(f,g)=
\begin{cases}
\mu^k(fg)+\Delta,
& \text{if } 0\leq k<M,\\[1ex]
\mu^k(fg)+\Delta+\dfrac{(k-M+2)(k-M+1)}{2},
& \text{if } M\leq k<m_0,\\[2ex]
\mu^k(fg)+\Delta+
\dfrac{(k-M+2)(k-M+1)-(k-m_0+1)(k-m_0)}{2},
& \text{if } k\geq m_0.
\end{cases}\]
\end{enumerate}

\end{remark}

The following example illustrates
Proposition~\ref{multfol} and Theorem~\ref{kbifu}.

\begin{example}\label{exbifu}\mbox{}
\begin{enumerate}
\item Let $f=x^2+y^2$ and
$g=x^2+2y^2+x^3+y^3$, so $m=\ell=2$.
The member $\alpha f+\beta g$ of the pencil has quadratic part
$(\alpha+\beta)x^2+(\alpha+2\beta)y^2$, so it is a node
unless $[\alpha:\beta]\in\{[1:-1],[2:-1]\}$.
The fibers of the pencil $\alpha f+\beta g$
are nodes, except
\[
h_1=f-g=-(y^2+x^3+y^3),\qquad
h_2=f-\tfrac12 g=\tfrac12(x^2-x^3-y^3),
\]
which are cusps. Hence
$\mathcal B^*(f,g)=\{[1:-1],[2:-1]\}$, $\mu(h_1)=\mu(h_2)=2$ and
$\mu_{\mathrm{gen}}=1$, so
$\Delta=(2-1)+(2-1)=2$.
Since $\nu(h_1)=\nu(h_2)=2=\nu_{\mathrm{gen}}$, we have $\delta=0$, and
Proposition~\ref{multfol} gives $m_0=m+\ell-1=3$. Indeed, the initial forms of $P=gf_x-fg_x$ and $Q=gf_y-fg_y$ are $2xy^2$ and $-2x^2y$, respectively.
This example shows that a bifurcation fiber may have the generic multiplicity.
Since $f$ and $g$ are nodes whose tangent cones have no common line,
$\mu(f)=\mu(g)=1$ and $I_0(f,g)=4$, so
$\mu(fg)=\mu(f)+\mu(g)+2I_0(f,g)-1=9$.
By Theorem~\ref{kbifu},
\[
\mu^k(f,g)=\mu^k(fg)+2
\quad\text{for all } k\geq 0,
\]
where

\[
\mu^k(fg)=9+k(k+1)-\varepsilon_k(3)
\quad\text{and}\quad
\mu^k(f,g)=11+k(k+1)-\varepsilon_k(3),
\]

\noindent with $\varepsilon_k(3)=0$ if $0\leq k<3$ and
$\varepsilon_k(3)=\frac{(k-2)(k-3)}{2}$ if $k\geq 3$. 
Since $\delta=0$, the two formulas have the
same correction term $\varepsilon_k(3)$, and
their difference is the constant $\Delta=2$.

\item Let $f=xy+x^3$ and $g=xy+y^3$,
so $m=\ell=2$. Here $fg=xy(y+x^2)(x+y^2)$, and
$\mu(fg)=\mu(f)+\mu(g)+2I_0(f,g)-1=13$,
since $\mu(f)=\mu(g)=1$ and $I_0(f,g)=6$.
The member $\alpha f+\beta g=(\alpha+\beta)xy+\alpha x^3+\beta y^3$
is a node unless $\alpha+\beta=0$.
The fibers of the pencil are nodes, except
$h=f-g=x^3-y^3$, three lines with
$\mu(h)=4$. Hence $\mathcal B^*(f,g)=\{[1:-1]\}$,
$\mu_{\mathrm{gen}}=1$ and
$\Delta=4-1=3$.
Since $\nu(h)=3$ and $\nu_{\mathrm{gen}}=2$, we have $\delta=1$, and
Proposition~\ref{multfol} gives $m_0=m+\ell=4$. Indeed, the initial forms of $P$ and $Q$ are $2x^3y+y^4$ and $-2xy^3-x^4$, respectively.
By Theorem~\ref{kbifu},
\[
\mu^k(f,g)=
\begin{cases}
\mu^k(fg)+3, & 0\leq k<4,\\[1ex]
\mu^k(fg)+k, & k\geq 4,
\end{cases}
\]
where

\[
\mu^k(fg)=\begin{cases}
13+k(k+1), & 0\le k<4,\\[1ex]
\dfrac{k^2+7k+20}{2}, & k\ge 4,
\end{cases}
\quad\text{and}\quad
\mu^k(f,g)=\begin{cases}
16+k(k+1), & 0\le k<4,\\[1ex]
\dfrac{k^2+9k+20}{2}, & k\ge 4.
\end{cases}
\]
Some values:

\begin{center}
\renewcommand{\arraystretch}{1.2}
\begin{tabular}{ccc}
\toprule
$k$ & $\mu^k(f,g)$ & $\mu^k(fg)$\\
\midrule
0 & 16 & 13\\
1 & 18 & 15\\
2 & 22 & 19\\
3 & 28 & 25\\
4 & 36 & 32\\
5 & 45 & 40\\
6 & 55 & 49\\
\bottomrule
\end{tabular}
\end{center}

\item Let $f=xy+x^3$ and
$g=xy+x^3+x^4+y^5$, so $m=\ell=2$.
Here
$\mu(fg)=\mu(f)+\mu(g)+2I_0(f,g)-1=19$,
since $\mu(f)=\mu(g)=1$ and $I_0(f,g)=9$.
The member $\alpha f+\beta g=(\alpha+\beta)(xy+x^3)+\beta(x^4+y^5)$
is a node unless $\alpha+\beta=0$. The fibers of the pencil
are nodes, except $h=f-g=-(x^4+y^5)$, with
$\mu(h)=12$. Hence
$\mathcal B^*(f,g)=\{[1:-1]\}$,
$\mu_{\mathrm{gen}}=1$ and
$\Delta=12-1=11$.
Since $\nu(h)=4$ and $\nu_{\mathrm{gen}}=2$, we have $\delta=2$, and
Proposition~\ref{multfol} gives $m_0=m+\ell+1=5$. Indeed, the initial forms of $P$ and $Q$ are $-3x^4y$ and $x^5$, respectively.
By Theorem~\ref{kbifu},
\[
\mu^k(f,g)=
\begin{cases}
\mu^k(fg)+11, & 0\leq k<4,\\[1ex]
\mu^k(fg)+12, & 4\leq k<5,\\[1ex]
\mu^k(fg)+2k+4, & k\geq 5,
\end{cases}
\]
where

\[
\mu^k(fg)=\begin{cases}
19+k(k+1), & 0\leq k<4,\\[1ex]
\dfrac{k^2+7k+32}{2}, & k\geq 4,
\end{cases}
\quad\text{and}\quad
\mu^k(f,g)=\begin{cases}
30+k(k+1), & 0\le k<5,\\[1ex]
\dfrac{k^2+11k+40}{2}, & k\ge 5.
\end{cases}
\]

\end{enumerate}
\end{example}

In \cite{Bodin1} and \cite{Sz2}, a
meromorphic germ
$(f:g):(\mathbb{C}^2,0)\dashrightarrow
\mathbb{P}^1$, with $f,g\in\mathfrak{m}$ relatively prime and $\mu(f,g)<\infty$, is called \emph{semitame}
if its only possible bifurcation fibers
are $f$ and $g$, that is,
$\mathcal{B}^*(f,g)=$\,$\emptyset$.
As a consequence of Proposition~\ref{multfol} and Theorem~\ref{kbifu},
the following corollary characterizes
semitame pairs in terms of the $k$-th
Milnor number.

\begin{corollary}\label{semitame}
Let $f,g\in\mathfrak{m}\subset\mathcal{O}_2$ be relatively
prime reduced germs
such that $\mu(f,g)<\infty$. Then
$\mu^k(f,g)\geq\mu^k(fg)$ for every integer $k\geq 0$, and the following
conditions are equivalent:
\begin{enumerate}
\item $(f:g)$ is semitame;
\item $\mu^k(f,g)=\mu^k(fg)$ for every integer $k\geq 0$;
\item $\mu^k(f,g)=\mu^k(fg)$ for some integer $k\geq 0$.
\end{enumerate}
\end{corollary}

\begin{proof}
For a fixed $k$, $\varepsilon_k(a)$ is a
non-increasing function of $a$.
Since $\delta\geq 0$, Theorem~\ref{kbifu}
gives
$\mu^k(f,g)-\mu^k(fg)=\Delta+
\varepsilon_k(M-1)-\varepsilon_k(M-1+\delta)
\geq\Delta\geq 0$
for every $k\geq 0$, where every summand of
$\Delta$ is positive. This proves the
inequality.

Assume first that $(f:g)$ is semitame, that
is, $\mathcal{B}^*(f,g)=\emptyset$. Then
$\Delta=\delta=0$, and the formula above gives
$\mu^k(f,g)=\mu^k(fg)$ for every $k\geq 0$,
which is the second condition. The second
condition clearly implies the third. Finally,
suppose that $\mu^k(f,g)=\mu^k(fg)$ for some
$k$. By the inequality above, $\Delta\leq 0$,
so $\Delta=0$. Since every summand of $\Delta$
is positive, $\mathcal{B}^*(f,g)=\emptyset$,
that is, $(f:g)$ is semitame. Hence the three
conditions are equivalent.
\end{proof}

\begin{remark}
By Proposition~\ref{multfol}, if $(f:g)$ is semitame then
$\nu(\mathcal{F}_{f/g})=\nu(f)+\nu(g)-1$; that is, the foliation
$\mathcal{F}_{f/g}$ has the smallest possible algebraic multiplicity.
The converse does not hold, as Example~\ref{exbifu}(1) shows.
\end{remark}

\begin{example}
Let $f=x^2+y^3$ and $g=x^3+y^2$. Here
$m=\ell=2$ and $m_0=3=m+\ell-1$, since the initial forms of $P=gf_x-fg_x$ and $Q=gf_y-fg_y$ are $2xy^2$ and $-2x^2y$, respectively.
The member $\alpha f+\beta g$ with $\alpha\beta\neq 0$ has quadratic part
$\alpha x^2+\beta y^2$, so it is a node. Hence
$\mu_{\mathrm{gen}}=1$, $\mathcal{B}^*(f,g)=\emptyset$, that is,
$(f:g)$ is semitame, and $\Delta=0$.
Since $f$ and $g$ are cusps with $I_0(f,g)=4$, we have
$\mu(f,g)=\mu(fg)=2+2+2\cdot 4-1=11$.

By Corollary~\ref{semitame},
\[
\mu^k(f,g)=\mu^k(fg)=
\begin{cases}
11+k(k+1), & 0\le k<3,\\[2mm]
11+k(k+1)-\dfrac{(k-2)(k-3)}{2},
& k\ge 3.
\end{cases}
\]
\end{example}

\vspace{1cm}
\noindent{\bf Acknowledgments:}
Arturo Fern\'andez-P\'erez and
Nancy Saravia-Molina would like to
thank the Universidad de La Laguna,
where this work was completed,
for its hospitality.
\\
\vspace{0.25cm}

\noindent{\bf Funding:}
The first author acknowledges support
from CNPq Projeto Universal
408687/2023-1 ``Geometria das
Equa\c{c}\~oes Diferenciais
Alg\'ebricas'' and from CNPq-Brazil
PQ-306011/2023-9.
The second author is supported by the
Spanish grant PID2023-149508NB-I00,
funded by MICIU/AEI/10.13039/501100011033
and by FEDER, UE.
The third author received financial
support from CAPES, Brazil.
The first and fourth authors were also
supported by the Vicerrectorado de
Investigaci\'on (VRI) of the Pontificia
Universidad Cat\'olica del Per\'u (PUCP)
through grant DFI-2025-PI1275.
\\
\vspace{0.25cm}

\noindent{\bf Data Availability Statement:}
Data sharing is not applicable to this
article, as no datasets were generated
or analyzed during the current study.
\vspace{0.25cm}
\\
\noindent{\bf Conflict of Interest:}
The authors declare that they have no
conflict of interest.

\end{document}